\documentclass{article}
\newsavebox{\mybox}
\usepackage{adjustbox}
\usepackage{subcaption}
\usepackage{comment}
\usepackage{amsmath}
\usepackage{listings}
\usepackage{mathrsfs,amsmath}
\usepackage{authblk}
\usepackage{indentfirst}
\usepackage{color,soul}
\usepackage{enumitem}
\usepackage{hyperref}
\usepackage{graphicx}
\usepackage{url}
\usepackage[section]{placeins}
\usepackage{appendix}
\usepackage{amssymb}
\usepackage{amsfonts}
\usepackage{listings}
\usepackage{color} 
\definecolor{mygreen}{RGB}{28,172,0} 
\definecolor{mylilas}{RGB}{170,55,241}
\usepackage{amstext}
\usepackage{amsthm}
\usepackage{appendix}
\DeclareMathOperator{\sech}{sech}

\title{On a series of Ramanujan, dilogarithm values, and solitons}
\author{Khristo Boyadzhiev\thanks{\textcolor{blue}{\underline{\href{mailto:k-boyadzhiev@onu.edu}{k-boyadzhiev@onu.edu}}}} \hspace{1mm}and Steven Manns\thanks{\textcolor{blue}{\underline{\href{mailto:s-manns.1@onu.edu}{s-manns.1@onu.edu}}}}}
\affil[]{Department of Mathematics}
\affil{Ohio Northern University, Ada, OH 45810, USA}
\date{\today}
\begin{document}

\maketitle

\begin{abstract} 
     We discuss several topics related to polylogarithms with focus on dilogarithms. The topics are: a generating function with harmonic numbers coming from Ramanujan, extending the dilogarithm to complex numbers beyond the unit disk, and relating polylogarithms to Bernoulli polynomials and to solitons. Several improvements and corrections are made to existing representations. 
     \vspace{5mm}
     
     \hspace{-\parindent}{\bf 2020 Mathematics subject classification:} 11M41, 33B30.\\
     \hspace{-\parindent}{\bf Key words:} Harmonic number, generating function, dilogarithm, trilogarithm, polylogarithm, Bernoulli polynomial, Riemann zeta function, solitons
\end{abstract}

\section{Introduction} \label{intro}
    \hspace{-\parindent}The polylogarithm function 
    \begin{eqnarray*}
        \operatorname{Li}_p(t) = \sum_{n=1}^\infty \frac{t^n}{n^p}
    \end{eqnarray*}
    has important applications in mathematics, especially in evaluating Euler type series and logarithmic integrals \cite{Bateman1953, Berndt1985, Kalman2012, Lewin1981, Mezo2014, Morris1979, Pickard1968, Srivastava2001, Srivastava2011, Srivastava2019, Stewart2020}. The function is continuous on the closed unit disk $|t| \leq 1$ and analytic in its interior. Especially popular are the dilogarithm and trilogarithm functions
    \begin{eqnarray*}
        \operatorname{Li}_2(t) = \sum_{n=1}^\infty \frac{t^n}{n^2}, \hspace{5mm} \operatorname{Li}_3 (t) = \sum_{n=1}^\infty \frac{t^n}{n^3}.
    \end{eqnarray*}
    Note that 
    \begin{eqnarray*}
        \operatorname{Li}_2(1) = \zeta(2) = \frac{\pi^2}{6}, \hspace{5mm} \operatorname{Li}_p(1) = \zeta(p)
    \end{eqnarray*}
    where $\zeta(s)$ is Riemann's zeta function (see \cite{Bateman1953, Berndt1985, Lewin1981, Srivastava2001, Srivastava2011, Srivastava2019, Zagier2007}). These functions appear in the structure of many generating functions involving the harmonic numbers \cite{Berndt1985, Dattoli2008, Hansen1975, Mezo2014, Srivastava2011, Stewart2020}
    \begin{eqnarray*}
        H_n = 1 + \frac{1}{2} + \cdots + \frac{1}{n}, \hspace{5mm} H_0 = 0.
    \end{eqnarray*}

    \hspace{-\parindent}The motivation for these notes came from the form of the generating function
    \begin{align}
        &\sum_{n=1}^\infty H_n \frac{t^{n+1}}{(n+1)^2} \nonumber\\
        & =  \frac{1}{2} \log (t) \log^2(1-t) + \log(1-t) \operatorname{Li}_2(1-t) - \operatorname{Li}_3(1-t) + \zeta (3)
       \label{gen_fun_raman}
    \end{align}
    where $0<t<1$. This generating function was computed by Ramanujan in his notebooks. It appears in Berndt \cite[p. 251]{Berndt1985}. In this form it also appears in Lewin's book \cite[p. 303]{Lewin1981}.\\
    
    \hspace{-\parindent}The right hand side in (\ref{gen_fun_raman}) is also defined for $t=0$ and for $t=1$ through the limits
    \begin{eqnarray*}
       \lim_{t \to 0^+} \log (t) \log^2(1-t) = 0 \hspace{5mm} \text{and} \hspace{5mm} \lim_{t \to 1^-} \log(1-t) \operatorname{Li}_2(1-t) = 0.
    \end{eqnarray*}
    Now over the interval $0 \leq t \leq 1$ the formula looks very good, for $t=1$ it readily gives the well-known value
    \begin{eqnarray*}
        \sum_{n=1}^\infty \frac{H_n}{(n+1)^2} = \zeta (3).
    \end{eqnarray*}
    However, this generating function has certain weaknesses. We cannot plug there numbers from the interval $-1 \leq t < 0$ without additional arrangements. At the same time, the left hand side is defined very well on the closed unit disk $|t| \leq 1$ being an absolutely convergent power series on this disk. Let us set
    \begin{eqnarray*}
        F(z)  = \sum_{n=0}^\infty H_n \frac{z^{n+1}}{(n+1)^2}.
    \end{eqnarray*}
    In order to compute its value for $z= - 1$, that is
    \begin{eqnarray*}
       F(-1) = \sum_{n=1}^\infty \frac{(-1)^{n+1} H_n}{(n+1)^2}
    \end{eqnarray*}
    we need to evaluate $\log(-1)$ and also $\operatorname{Li}_2(2),$  $\operatorname{Li}_3(2)$. We will discuss these values later, now we return to the series in (\ref{gen_fun_raman}). How do we come to this representation? \\
    
    \hspace{-\parindent}The series (\ref{gen_fun_raman}) comes from the well-known logarithmic series
    \begin{eqnarray*}
       -\log(1-t) = \sum_{n=1}^\infty \frac{t^n}{n}, \hspace{5mm} - 1 \leq t < 1
    \end{eqnarray*}
    which we multiply by the geometric series
    \begin{eqnarray*}
        \frac{1}{1-t} = \sum_{n=0}^\infty t^n, \hspace{5mm} |t|<1
    \end{eqnarray*}
    to get the well-known generating function for the harmonic numbers
    \begin{eqnarray*}
        - \frac{\log(1-t)}{1-t} = \sum_{n=0}^\infty H_n t^n, \hspace{5mm} |t|<1.
    \end{eqnarray*}
    Integration gives
    \begin{eqnarray*}
        \sum_{n=0}^\infty H_n \frac{t^{n+1}}{n+1} = \frac{1}{2} \log^2(1-t), \hspace{5mm} \text{or} \hspace{5mm} \sum_{n=0}^\infty H_n \frac{t^n}{n+1} = \frac{\log^2(1-t)}{2t}
    \end{eqnarray*}
    and from this (introducing the complex variable $z$)
    \begin{eqnarray}
       \sum_{n=0}^\infty H_n \frac{z^{n+1}}{(n+1)^2} = \int_0^z \frac{\log^2(1-x)}{2x} dx = F(z).
       \label{func_F}
    \end{eqnarray}
    The integral here represents the analytic function $F(z)$ on the unit disk, continuous on the boundary. The left hand side is its Taylor series. For $0<t<1$ the integral has the explicit representation given in (\ref{gen_fun_raman}). For $-1<t<1$ we have
    \begin{eqnarray*}
       F(-t)\hspace{-2mm} & = & \hspace{-3mm}\int_0^{-t} \frac{\log^2(1-x)}{2x} dx = \int_0^t \frac{\log^2(1+x)}{2x} dx \\ \\
       & = &\hspace{-2mm} \frac{1}{2} \log(t) \log^2(1+t) - \frac{1}{3} \log^3(1+t) - \log(1+t) \operatorname{Li}_2 \biggr{(} \frac{1}{1+t} \biggr{)} \\ \\
       &&- \operatorname{Li}_3 \biggr{(} \frac{1}{1+t} \biggr{)} + \zeta(3)
    \end{eqnarray*}
    as given by Lewin \cite[p. 310]{Lewin1981}. This way we have for $0< t \leq 1$
    \begin{align}
        &\sum_{n=0}^\infty (-1)^{n+1} H_n\frac{t^{n+1}}{(n+1)^2} \nonumber\\ \nonumber\\
        &= \frac{\log(t) \log^2(1+t)}{2} - \frac{\log^3(1+t)}{3} - \log(1+t) \operatorname{Li}_2 \biggr{(} \frac{1}{1+t} \biggr{)}- \operatorname{Li}_3 \biggr{(} \frac{1}{1+t} \biggr{)} + \zeta(3).
        \label{f_at_minus1}
    \end{align}
    We can now plug here $t=1$ to get a well-defined number
    \begin{eqnarray*}
        \sum_{n=1}^\infty \frac{(-1)^{n+1} H_n}{(n+1)^2} = - \frac{1}{3}(\log2)^3 - \log (2) \operatorname{Li}_2 \biggr{(} \frac{1}{2} \biggr{)} - \operatorname{Li}_3 \biggr{(} \frac{1}{2} \biggr{)} + \zeta(3)
    \end{eqnarray*}
    (this will be specified below). It turns out that the analytic function $F(z)$ defined in (\ref{func_F}) has two faces like the Roman god Janus - one for positive and one for negative numbers. This curious fact is a good reason to study the extensions of $\operatorname{Li}_2(x)$ and $\operatorname{Li}_3(x)$ beyond the unit disk and find various explicit values. First of all, it is desirable to express $F(t)$ by a single formula for the entire interval $-1 \leq t \leq 1$. This problem is addressed in the next section.\\
    
    \hspace{-\parindent}In the third section we discuss particular values of $\operatorname{Li}_2(t)$, $\operatorname{Li}_3 (t)$, extensions beyond the unit disk, and give a slight improvement of Lewin's formula for the complex dilogarithm. In the fourth section we collect some observations on the values of $\operatorname{Li}_2 (t)$ and discuss open problems. In the fifth and last section we give a new proof of the important formula \cite[(7.192)]{Lewin1981} 
    \begin{eqnarray*}
        \operatorname{Li}_p(z) + (-1)^p \operatorname{Li}_{p} \biggr{(} \frac{1}{z} \biggr{)} = A_p B_p \biggr{(} \frac{\log(z)}{2 \pi i} \biggr{)}
    \end{eqnarray*}
    and relate it to the soliton solution $u(x,t) = -2 \sech^2(x-4t)$ of the Korteweg-de Vries equation $u_t - 6u u_x + u_{xxx} = 0$. Here $B_p(x)$ are the Bernoulli polynomials and $A_p$ is a constant.

\section{A new generating function}\label{new_generating}
    \hspace{-\parindent}We want to modify the generating function in (\ref{gen_fun_raman}) so that it could be evaluated for all $-1 \leq t \leq 1$. For this purpose we need to remove the first term and also $\operatorname{Li}_2(1-t)$ and $\operatorname{Li}_3(1-t)$. In order to do this we will use two important identities. The first one is the reflection formula
    \begin{eqnarray}
       \operatorname{Li}_2 (1-t) = \frac{\pi^2}{6} - \operatorname{Li}_2(t) - \log(t) \log(1-t) \hspace{5mm} (0 \leq t \leq 1)
       \label{reflection_form}
    \end{eqnarray}
    which comes from Euler - see comments in \cite{Kalman2012}. The second one is Landen's identity for the trilogarithm \cite[p. 296]{Lewin1981}, \cite[p. 113]{Sitaramachandrarao1987}, \cite[p.181]{Srivastava2001}
    \begin{align}
        &\operatorname{Li}_3 \biggr{(} \frac{-t}{1-t}\biggr{)} \nonumber\\
        &= \frac{1}{6} \log^3(1-t) - \operatorname{Li}_3(1-t) - \frac{1}{2} \log(t) \log^2(1-t) + \frac{\pi^2}{6} \log(1-t) - \operatorname{Li}_3(t) + \zeta(3).
        \label{landen_trilog}
    \end{align}
    Before proceeding further we want to make one observation. From (\ref{landen_trilog}) we find
    \begin{eqnarray*}
       \lim_{t \to 1^-} \Biggr{[}\operatorname{Li}_3 \biggr{(} \frac{-t}{1-t} \biggr{)} - \frac{1}{6} \log^3(1-t) - \frac{\pi^2}{6} \log(1-t)  \Biggr{]} = 0
    \end{eqnarray*}
    and from here also
    \begin{eqnarray*}
        \lim_{t \to 1^-} \Biggr{[}\operatorname{Li}_3 \biggr{(} \frac{-t}{1-t} \biggr{)} - \frac{1}{6} \log^3(1-t) - \operatorname{Li}_2(t) \log(1-t) \Biggr{]} = 0
    \end{eqnarray*}
    since $\frac{\pi^2}{6} = \operatorname{Li}_2(1)$ and 
    \begin{eqnarray*}
        \lim_{t \to 1^-} \biggr{[}\biggr{(}\operatorname{Li}_2(t) - \operatorname{Li}_2(1) \biggr{)}\log(1-t) \biggr{]}= 0
    \end{eqnarray*}
    (a simple calculus exercise).\\
    
    \hspace{-\parindent}Thus we come to the following proposition.\\
    
    \hspace{-\parindent}\textbf{Proposition 1.} For every $-1 \leq t \leq 1$
    \begin{align}
        F(t) &=\sum_{n=1}^\infty H_n \frac{t^{n+1}}{(n+1)^2} \nonumber\\ \nonumber \\
       &=  \operatorname{Li}_3 \biggr{(} \frac{-t}{1-t}\biggr{)} - \frac{1}{6} \log^3(1-t) - \log(1-t) \operatorname{Li}_2(t) + \operatorname{Li}_3(t).
       \label{prop1}
    \end{align}
    \textit{Proof of the proposition.} First, it is easy to see that $-1 \leq t \leq \frac{1}{2}$ if and only if $-1 \leq \frac{-t}{1-t} \leq \frac{1}{2}$, so that all terms on the right hand side in (\ref{prop1}) are well defined on this interval. For $\frac{1}{2}\leq t < 1$ we have $-\infty < \frac{-t}{1-t} \leq -1$ and for such $t$ we use Landen's identity (\ref{landen_trilog}) to evaluate $\operatorname{Li}_3 (\frac{-t}{1-t})$. The value at $t=1$ (as a limit) is $\zeta(3)$ in view of the above observation. Landen's identity practically shows that the right hand side in (\ref{prop1}) is well defined for all $-1 \leq t \leq 1$.\\
    
    \hspace{-\parindent}In the representation (\ref{gen_fun_raman}) we first replace the term $\operatorname{Li}_2(1-t)$ by using the reflection formula for the dilogarithm (\ref{reflection_form}) 
    \begin{align}
        &F(t)= \frac{1}{2} \log(t) \log^2(1-t) + \log(1-t) \biggr{[}\frac{\pi^2}{6} - \operatorname{Li}_2(t) - \log (t) \log(1-t)\biggr{]} \nonumber\\
        &- \operatorname{Li}_3(1-t) + \zeta(3).
    \end{align}
    That is,
    \begin{eqnarray*}
        F(t) = - \frac{1}{2} \log(t) \log^2(1-t) + \frac{\pi^2}{6} \log(1-t) - \log(1-t) \operatorname{Li}_2(t) - \operatorname{Li}_3(1-t) + \zeta(3).
    \end{eqnarray*}
    Next we use Landen's identity (\ref{landen_trilog}) to express $\operatorname{Li}_3(1-t)$ from there
        \begin{align}
            & \operatorname{Li}_3(1-t) = \frac{1}{6} \log^3(1-t) - \operatorname{Li}_3 \biggr{(} \frac{-t}{1-t} \biggr{)} - \frac{1}{2} \log(t) \log^2(1-t) + \frac{\pi^2}{6} \log(1-t) \nonumber\\  
            &- \operatorname{Li}_3(t) + \zeta(3) 
            \label{eqn7_draft}
        \end{align}
    and replace it in the above representation of $F(t)$. This gives the desired equation (\ref{prop1}). As an added bonus we got rid of the term with $\log(t) \log^2(1-t)$.\\
    
    \hspace*{\fill} $\blacksquare$

\section{Extending the dilogarithm and the trilogarithm beyond the unit disk} \label{extending_dilog}
    \hspace{-\parindent}First we list some well-known values
    \begin{eqnarray}
       \operatorname{Li}_2(-1) = - \frac{1}{2} \zeta(2) = - \frac{\pi^2}{12}; \hspace{5mm} \operatorname{Li}_3(-1) = - \frac{3}{4} \zeta(3)
       \label{dilog_trilog_atneg1}
    \end{eqnarray}
    \begin{eqnarray}
        \operatorname{Li}_2 \biggr{(}\frac{1}{2} \biggr{)} = \frac{\pi^2}{12} - \frac{(\log 2)^2}{2} \label{dilog_at_half}
    \end{eqnarray}
    (which comes immediately from the reflection formula (\ref{reflection_form}) with $t = \frac{1}{2}$)
    \begin{eqnarray}
       \operatorname{Li}_3 \biggr{(} \frac{1}{2}\biggr{)} = \frac{7}{8} \zeta(3) - \frac{\pi^2}{12} \log (2) + \frac{1}{6}(\log 2)^3
       \label{trilog_at_half}
    \end{eqnarray}
    (resulting from (\ref{eqn7_draft}) with $t = \frac{1}{2}$ and (\ref{dilog_trilog_atneg1})).\\
    
    \hspace{-\parindent}Using these numbers we can evaluate the series in (\ref{gen_fun_raman}) for $t = \frac{1}{2}$. We can use either the original representation (\ref{gen_fun_raman}) or (\ref{prop1}) to compute
    \begin{eqnarray}
       \sum_{n=1}^\infty \frac{H_n}{2^{n+1}(n+1)^2} = \frac{1}{8} \zeta(3) - \frac{1}{6} (\log 2)^3.
       \label{eqn11_draft}
    \end{eqnarray}
    Also from (\ref{f_at_minus1})
    \begin{eqnarray}
        \sum_{n=1}^\infty \frac{(-1)^{n+1} H_n}{(n+1)^2} = \frac{1}{8} \zeta(3).
        \label{eqn12_draft}
    \end{eqnarray}
    From here using the equation
    \begin{eqnarray*}
        \sum_{n=1}^\infty \frac{(-1)^{n+1} H_n}{(n+1)^2} = \sum_{n=0}^\infty \biggr{(}H_{n+1} - \frac{1}{n+1} \biggr{)} \frac{(-1)^{n+1}}{(n+1)^2} = \sum_{n=1}^\infty \frac{(-1)^n H_n}{n^2} - \operatorname{Li}_3(-1)
    \end{eqnarray*}
    we find also
    \begin{eqnarray}
       \sum_{n=1}^\infty \frac{(-1)^{n-1} H_n}{n^2} = \frac{5}{8} \zeta(3).
       \label{eqn13_draft}
    \end{eqnarray}
    This evaluation was published by Sitaramachandrarao \cite{Sitaramachandrarao1987} in 1987.\\
    
    \hspace{-\parindent}However, we cannot evaluate explicitly the series (\ref{gen_fun_raman}) for $t = - \frac{1}{2}$ at this point, because in (\ref{f_at_minus1}) this brings to the values $\operatorname{Li}_2 (\frac{2}{3})$ and $\operatorname{Li}_3(\frac{2}{3})$. Using (\ref{prop1}) for $t = - \frac{1}{2}$ we come also to $\operatorname{Li}_2(-\frac{1}{2})$ and $\operatorname{Li}_3(-\frac{1}{2})$.\\
    
    \hspace{-\parindent}The extension of the dilogarithm beyond the unit disk was considered by Lewin \cite[Chapter 5]{Lewin1981}. We will briefly review some cases here. In the Introduction we mentioned the values $\operatorname{Li}_2(2)$ and $\operatorname{Li}_3(2)$ and we will review them first. To define $\operatorname{Li}_2 (2)$ we start with the reflection formula (\ref{reflection_form})
    \begin{eqnarray*}
       \operatorname{Li}_2(1-x) = \frac{\pi^2}{6} - \operatorname{Li}_2(x) - \log(x) \log(1-x) \hspace{5mm} (0 \leq 
       x \leq 1).
    \end{eqnarray*}
    The last term approaches zero when $x \to 0^+$ and when $x \to 1^-$, so it is considered zero for $x=0$ and for $x=1$. We see that the extension of $\operatorname{Li}_2(1-x)$ depends on the complex logarithm $\log (x)$ and this extension, in general, will be multivalued. We will use only the principal branch of the logarithm. For any $z \in \mathbb{C} - (-\infty, 0]$ we have
    \begin{eqnarray}
       \log (z)  = \ln |z| + i \operatorname{Arg}(z), \hspace{5mm} - \pi < \operatorname{Arg}(z) \leq \pi
       \label{complex_log}
    \end{eqnarray}
    where $\operatorname{Arg}(z)$ is the principal value of the argument. The numbers on the negative half-axis $(-\infty, 0]$ are written as $r e^{i \pi}$, $r>0$. The principal branch of the logarithm on $(-\infty,0]$ has the form
    \begin{eqnarray}
       \log(r e^{i\pi}) = \ln (r) + i \pi, \hspace{5mm} r>0
       \label{principal_log_neg_reals}
    \end{eqnarray}
    so that $\log(-1) = i \pi$. This will be fixed throughout the paper.\\
    
    \hspace{-\parindent}Now in the reflection formula (\ref{reflection_form}) we can move $x$ through the origin into the interval $[-1,0]$. For $x = r e^{i \pi}$ $(0 \leq r \leq 1)$ we have
    \begin{eqnarray}
       \operatorname{Li}_2(1-r e^{i\pi}) = \frac{\pi^2}{6} - \operatorname{Li}_2(-r) - (\log (r) + i \pi) \log(1+r)
       \label{eqn15_draft}
    \end{eqnarray}
    and with $r \to 1^-$ we find
    \begin{eqnarray*}
        \operatorname{Li}_2(2) = \frac{\pi^2}{6} - \operatorname{Li}_2(-1) - i \pi \log (2) = \frac{\pi^2}{4} - i \pi \log (2)
    \end{eqnarray*}
    that is,
    \begin{eqnarray}
       \operatorname{Li}_2(2) = \frac{\pi^2}{4} - i \pi \log (2) .
       \label{dilog_at2}
    \end{eqnarray}
    This value of $\operatorname{Li}_2(2)$ has been adopted by the mathematical community. It is confirmed by another identity from \cite[p.283]{Lewin1981}.
    \begin{eqnarray}
       \operatorname{Li}_2(x) + \operatorname{Li}_2 \biggr{(} \frac{1}{x} \biggr{)} = \frac{\pi^2}{3} - \frac{\ln^2 (x)}{2} - i \pi \log (x) \hspace{5mm} (1<x)
       \label{eqn17_draft}
    \end{eqnarray}
    where for $x = 2$ we come to (\ref{dilog_at2}) again. This value is also cited by Morris \cite{Morris1979}.\\
    
    \hspace{-\parindent}In a similar manner we can approach the value $\operatorname{Li}_3(2)$ by using Landen's identity in the form (\ref{eqn7_draft}). Moving $t$ through the origin toward the point $-1 = e^{i\pi}$ and assuming that $\log(t) = \log|t| + i \pi$ and $\log(-1) = i \pi$ we come to the evaluation
    \begin{eqnarray*}
       \operatorname{Li}_3(2) = \frac{1}{6} \log^3 (2) - \operatorname{Li}_3 \biggr{(}\frac{1}{2} \biggr{)} - \frac{i \pi}{2} \log^2 (2) + \frac{\pi^2}{6} \log (2) - \operatorname{Li}_3(-1) + \zeta(3)
    \end{eqnarray*}
    which simplified is 
    \begin{eqnarray}
       \operatorname{Li}_3(2) = \frac{\pi^2}{4}\log (2) + \frac{7}{8} \zeta(3) - \frac{i \pi}{2} \log^2 (2).
       \label{eqn18_draft}
    \end{eqnarray}
    Lewin \cite[p.296]{Lewin1981} lists the identity
    \begin{eqnarray}
       \operatorname{Li}_3(x) = \operatorname{Li}_3 \biggr{(}\frac{1}{x} \biggr{)} + \frac{\pi^2}{3} \log (x) - \frac{1}{6} \log^3 (x) - \frac{i \pi}{2} \log^2 (x) \hspace{5mm} (1<x)
       \label{eqn19_draft}
    \end{eqnarray}
    which for $x=2$ implies the same result. This identity provides an extension of $\operatorname{Li}_3(x)$ for the entire half-plane $\operatorname{Re} x >1$.\\
    
    \hspace{-\parindent}The equation
    \begin{eqnarray*}
        \operatorname{Li}_2 \biggr{(} \frac{-x}{1-x} \biggr{)} = - \frac{1}{2} \log^2(1-x) - \operatorname{Li}_2(x)
    \end{eqnarray*}
    (easily verified by differentiation) is Landen's identity for the dilogarithm \cite[p.5]{Lewin1981}. It can be written in the form
    \begin{eqnarray}
       \operatorname{Li}_2(-z) = - \operatorname{Li}_2 \biggr{(} \frac{z}{1+z} \biggr{)} - \frac{1}{2} \log^2 (1+z)
       \label{landen_dilog}
    \end{eqnarray}
    which conveniently provides an extension of the dilogarithm to the left half-plane $\operatorname{Re} z \leq 0$, because if $\operatorname{Re} z \geq 0$, then $|\frac{z}{1+z}| \leq 1$ and the two terms on the right hand side in (\ref{landen_dilog}) are well defined. This way
    \begin{eqnarray}
       \operatorname{Li}_2(-2) = - \operatorname{Li}_2\biggr{(} \frac{2}{3} \biggr{)} - \frac{1}{2} \log^2(3), \hspace{5mm} \operatorname{Li}_2(-3) = - \operatorname{Li}_2 \biggr{(} \frac{3}{4} \biggr{)} - \frac{1}{2} \log^2(4).
       \label{eqn21_draft}
    \end{eqnarray}
    Equation (\ref{landen_dilog}) gives also
    \begin{eqnarray}
       \operatorname{Li}_2 \biggr{(} - \frac{1}{2} \biggr{)} = - \operatorname{Li}_2 \biggr{(} \frac{1}{3} \biggr{)} - \frac{1}{2} \log^2 \biggr{(} \frac{3}{2}\biggr{)}.
    \end{eqnarray}
    Evaluations of these in the spirit of (\ref{dilog_trilog_atneg1}) and (\ref{dilog_at_half}) are presently unknown. \\
    
    \hspace{-\parindent}The dilogarithm can be defined also by the line integral
    \begin{eqnarray*}
       \operatorname{Li}_2(z) = - \int_0^z \frac{\log(1-t)}{t} dt
    \end{eqnarray*}
    or, which is better, by the ordinary integral
    \begin{eqnarray}
       \operatorname{Li}_2(-z) = - \int_0^1 \frac{\log(1+zt)}{t} dt.
       \label{dilog_integral}
    \end{eqnarray}
    For $|z|<1$ the representation (\ref{dilog_integral}) can be proved directly by expanding the logarithm and integrating term by term. Using the principal branch of the complex logarithm the representation (\ref{dilog_integral}) provides an extension of $\operatorname{Li}_2(-z)$ for all $z$ in the complex plane $\mathbb{C}$ cut along $(- \infty, -1]$.\\
    
    \hspace{-\parindent}Lewin \cite[p.121]{Lewin1981} gives the representation
    \begin{eqnarray}
       \operatorname{Li}_2(re^{i\theta}) = - \frac{1}{2}\int_0^r \frac{\log(1 - 2t \cos \theta + t^2)}{t} dt + i \int_0^r \arctan \biggr{(} \frac{y \sin \theta}{1 - y \cos \theta}\biggr{)} \frac{dy}{y}
       \label{lewin_real_imaginary}
    \end{eqnarray}
    for $z = r e^{i\theta}$ by using the real and imaginary parts of the complex logarithm. The same representation is used also by Osacar et al. in \cite{Osacar1995}. However, this formula is incomplete, because the range of the arctangent is only $(-\frac{\pi}{2}, \frac{\pi}{2})$ and so (\ref{lewin_real_imaginary}) represents only part of the principal argument of the complex logarithm. As noted by Baker and Sluis in \cite{baker1997}, for a complex number $z = x + iy$ the correct principal argument is given by
    \begin{eqnarray}
       \operatorname{Arg} (x+iy) = 
       \begin{cases}
           2 \arctan \biggr{(} \frac{y}{x + \sqrt{x^2 + y^2}} \biggr{)} & (x>0 \hspace{2mm} \text{or} \hspace{2mm} y =0) \\
           \pi & (x<0 \hspace{2mm}\text{and} \hspace{2mm}y = 0)
       \end{cases}
       \label{pricipal_arg_form}
    \end{eqnarray}
    with range $- \pi < \operatorname{Arg}(z) \leq \pi$. Now we state the following proposition based on (\ref{dilog_integral}) and (\ref{pricipal_arg_form}) which improves the imaginary part in (\ref{lewin_real_imaginary}). \\
    
    \hspace{-\parindent}\textbf{Proposition 2.} For all $z = x+iy = re^{i\theta}$ in the complex plane $\mathbb{C}$ cut along $(- \infty, -1]$ we have
    \begin{eqnarray}
       \operatorname{Li}_2(-z) \hspace{-2mm}&=&\hspace{-2mm} - \int_0^1 \frac{\ln\biggr{(}1+2xt + t^2(x^2+y^2)\biggr{)}}{2t} dt \nonumber \\ \nonumber\\
       &&\hspace{-2mm}- i \int_0^1 \frac{2}{t} \arctan \biggr{(} \frac{yt}{1 + xt + \sqrt{1 + 2xt + t^2(x^2+y^2)}}\biggr{)}dt
       \label{prop2_rep1}
    \end{eqnarray}
    or,
    \begin{eqnarray}
       \operatorname{Li}_2(-z) \hspace{-2mm}& = &\hspace{-2mm} - \int_0^1 \frac{\ln(1 + 2rt \cos \theta + t^2 r^2)}{2t} dt \nonumber \\ \nonumber\\
       &&\hspace{-2mm}- i \int_0^1 \frac{2}{t} \arctan \biggr{(} \frac{rt \sin \theta}{1 + rt \cos \theta + \sqrt{1 + 2rt \cos \theta + t^2 r^2}} \biggr{)} dt.
       \label{prop2_rep2}
    \end{eqnarray}
    With $x= 0$ and $|y| \leq 1$ we find from (\ref{prop2_rep1}) by expanding the logarithm
    \begin{eqnarray*}
       \operatorname{Re} \operatorname{Li}_2(-iy)  = - \int_0^1 \frac{\ln (1+t^2y^2)}{2t} dt = \frac{1}{4} \sum_{n=1}^\infty \frac{(-1)^n y^{2n}}{n^2} = \frac{1}{4} \operatorname{Li}_2(-y^2)
    \end{eqnarray*}
    confirming the result of Lewin \cite[p.38]{Lewin1981} on the values of the dilogarithm on the imaginary axis. Also from (\ref{prop2_rep1}) for $x=0$ and integrating by parts twice
    \begin{align*}
        \operatorname{Im} \operatorname{Li}_2(-iy) & =  -2 \int_0^1 \frac{1}{t} \arctan \biggr{(} \frac{yt}{1 + \sqrt{1 + t^2y^2}} \biggr{)} dt\\ \\
       & = - 2 \int_0^1 \arctan \biggr{(} \frac{yt}{1 + \sqrt{1 + t^2 y^2}} \biggr{)} d \log t \\ \\
       & =  y \int_0^1 \frac{\log (t)}{1 + t^2y^2} dt\\ \\
       & = - \int_0^1 \frac{\arctan(yt)}{t} dt; \\ \\
       \operatorname{Im} \operatorname{Li}_2(iy)& = \int_0^1 \frac{\arctan(yt)}{t} dt
    \end{align*}
    confirming equation (2.1) on p. 38 in \cite{Lewin1981}. In particular,
    \begin{eqnarray}
       \operatorname{Im} \operatorname{Li}_2(i) = G, \hspace{5mm} \operatorname{Im} \operatorname{Li}_2(-i) = -G
    \end{eqnarray}
    ($G$ - the Catalan constant). In the above computation $y$ is an arbitrary real number, while the equation (2.1) in \cite[p.38]{Lewin1981} requires $|y| \leq 1$.\\
    
    \hspace{-\parindent}In a second example we consider values on the line $y=x$, that is, $z = x+ix$. Lewin \cite[p.132]{Lewin1981} evaluated $\operatorname{Re} \operatorname{Li}_2 (-z)$. Here we will look at $\operatorname{Im} \operatorname{Li}_2(-z)$. Writing $\frac{1}{t} dt = d \log(t)$ and integrating by parts we put $\operatorname{Im} \operatorname{Li}_2 (-x-ix)$ in different forms
    \begin{eqnarray*}
        \operatorname{Im} \operatorname{Li}_2 (-x-ix) \hspace{-2mm}& = & \hspace{-2mm}- \int_0^1 \frac{2}{t} \arctan \biggr{(} \frac{xt}{1 + xt + \sqrt{1 + 2xt + 2x^2t^2}}\biggr{)} dt \\ \\
        & = &\hspace{-2mm} x\int_0^1 \frac{\log (t)}{1 + 2xt + 2x^2 t^2} dt = 2x \int_0^1 \frac{\log (t)}{(2xt+1)^2 + 1}dt \\ \\
        & = &\hspace{-2mm} \int_0^1 \log (t) d\biggr{(} \arctan(2xt+1) - \frac{\pi}{4} \biggr{)}\\ \\
        & = &\hspace{-2mm} \int_0^1 \biggr{(} \frac{\pi}{4} - \arctan(2xt+1) \biggr{)} \frac{dt}{t}.
    \end{eqnarray*}
    On the line $y = -x$ we have
    \begin{eqnarray*}
        \operatorname{Im} \operatorname{Li}_2(-x+ix) \hspace{-2mm}& = &\hspace{-2mm} - \operatorname{Im} \operatorname{Li}_2(-x-ix) \\ \\
        & = & \hspace{-2mm}- x \int_0^1 \frac{\log (t)}{1 + 2xt + 2x^2 t^2}dt  = \int_0^1 \biggr{(} \arctan(2xt+1) - \frac{\pi}{4} \biggr{)} \frac{dt}{t}.
    \end{eqnarray*}
    For related results see also \cite[p. 106-119]{Srivastava2001} and \cite[p. 175-190]{Srivastava2011}.    

\section{Some observations and open problems}\label{open_problems}
    \hspace{-\parindent}In a similar spirit to the work done for the dilogarithm in section \ref{extending_dilog}, we note the trilogarithm can be written in terms of its real and imaginary parts. On page 153 in Lewin's book \cite{Lewin1981}, the trilogarithm is defined by the integral
    \begin{eqnarray*}
       \operatorname{Li}_3(z) = \int_0^z \frac{\operatorname{Li}_2(z)}{z} dz,
    \end{eqnarray*}
    where this definition holds irrespective of whether $z$ is on the interior of the unit disk. However, an alternative form is given by
    \begin{eqnarray}
        \operatorname{Li}_3(z) = \int_0^1 \frac{\operatorname{Li}_2(zx)}{x}dx,
        \label{trilog_integral}
    \end{eqnarray}
    which can be proven for $|zx| \leq 1$ by using the series for $\operatorname{Li}_2(zx)$ and integrating term by term. Consider that the representation of the dilogarithm using (\ref{dilog_integral}) can be substituted into this integral. After making this substitution for the integral form of $\operatorname{Li}_3(-z)$ coming from (\ref{trilog_integral}), the result is a representation of $\operatorname{Li}_3(-z)$ as a double integral. More specifically, upon writing the complex logarithm in terms of its real and imaginary parts, we find for all $z = u + i v = r e^{i\theta}$
     \begin{align}\label{trilog_real_imag1}
        &\operatorname{Li}_3(-z)  \nonumber\\
        &= - \int_0^1 \int_0^1 \frac{\log \biggr{(} 1 + 2uxt + (u^2 + v^2)x^2t^2 \biggr{)}}{2xt}  dt dx \nonumber\\ \nonumber \\ 
        & - i \int_0^1 \int_0^1 \frac{2}{xt} \arctan \biggr{(} \frac{vxt}{1 + uxt + \sqrt{1 + 2uxt + (u^2 + v^2)x^2 t^2}} \biggr{)} dt dx 
    \end{align}
    or,
    \begin{align}
        &\operatorname{Li}_3(-z) \nonumber \\ 
        &=  - \int_0^1 \int_0^1 \frac{\log\biggr{(}1+2rxt \cdot\cos(\theta) + r^2x^2t^2\biggr{)}}{2xt}  dt dx \nonumber\\ \nonumber \\ \nonumber
        & - i \int_0^1 \int_0^1 \frac{2}{xt} \arctan \biggr{(} \frac{rxt \cdot \sin(\theta)}{1 + rxt \cdot \cos(\theta) + \sqrt{1+2rxt\cdot \cos(\theta) + r^2 x^2 t^2}} \biggr{)} dt dx .\\
       \label{trilog_real_imag2}
    \end{align}
    
    \hspace{-\parindent}By using the trigonometric identity
    \begin{eqnarray*}
       \arctan (x) = 2 \arctan\biggr{(} \frac{x}{1 + \sqrt{1 + x^2}}\biggr{)},
    \end{eqnarray*}
    which follows from the tangent half-angle formula, letting $z = - i$ in formula (\ref{trilog_real_imag1}) shows that
    \begin{eqnarray*}
       \operatorname{Li}_3(i) = - \int_0^1 \int_0^1 \frac{\log(1+x^2 t^2)}{2xt} dt dx + i \int_0^1 \int_0^1 \frac{\arctan(xt)}{xt} dt dx.
    \end{eqnarray*}
    Notice that $0 \leq x \leq 1$ and $0 \leq t \leq 1$, so both the logarithm and arctangent can be expanded using their Maclaurin series. After putting these series into the integrals and integrating term by term first with respect to $t$ and then with respect to $x$, the result is
    \begin{eqnarray*}
       \operatorname{Li}_3(i) = \frac{1}{8} \sum_{n=1}^\infty \frac{(-1)^n}{n^3} + i \sum_{n=0}^\infty \frac{(-1)^n}{(2n+1)^3}
    \end{eqnarray*}
    or,
    \begin{eqnarray}
       \operatorname{Li}_3(i) = - \frac{3}{32} \zeta(3) + \frac{3\pi i}{16} \zeta(2).
    \end{eqnarray}
    This value of $\operatorname{Li}_3(i)$ matches what is given in the Appendix of \cite{Stewart2020} by Stewart, which is a good sign our representation of the trilogarithm is indeed useful. \\
    
    \hspace{-\parindent}Turning now to problems related to these notes that remain open, one unknown value of the dilogarithm has caught our attention in particular. That is, given that $\operatorname{Li}_2(1)$, $\operatorname{Li}_2(-1)$, and $\operatorname{Li}_2(\frac{1}{2})$ are all known values ($\frac{\pi^2}{6}$, $-\frac{\pi^2}{12}$, and $\frac{\pi^2}{12} - \frac{\log^2(2)}{2}$, respectively), one may expect that $\operatorname{Li}_2(-\frac{1}{2})$ could be evaluated. However, in spite of the fact that obtaining a value for $\operatorname{Li}_2(-\frac{1}{2})$ would make many other values of the dilogarithm known, $\operatorname{Li}_2(-\frac{1}{2})$ has yet to be computed.  \\
    
    \hspace{-\parindent}Due to this observation, we would like to provide a source of motivation for others to seek to determine $\operatorname{Li}_2(-\frac{1}{2})$. Referring to this value from now on as the dilogarithm constant, let it be denoted by 
    \begin{eqnarray}
       d_2 = \operatorname{Li}_2 \biggr{(}-\frac{1}{2}\biggr{)}.
       \label{dilog_constant}
    \end{eqnarray}
    Setting $x= \frac{1}{2}$ in formula (13) of Lewin's book \cite[p.283]{Lewin1981} and using the known value for $\operatorname{Li}_2(\frac{1}{2})$, it is readily seen that 
    \begin{eqnarray}
       \operatorname{Li}_2\biggr{(} \frac{1}{4}\biggr{)} = 2 d_2 + \frac{\pi^2}{6} - \log^2(2).
       \label{dilog_one_fourth}
    \end{eqnarray}
    Likewise, taking $x=2$ in formula (5) of \cite[p.283]{Lewin1981} relates $d_2$ to $\operatorname{Li}_2(-2)$ by the equation
    \begin{eqnarray}
       \operatorname{Li}_2 (-2) = - d_2 - \frac{\pi^2}{6} - \frac{1}{2} \log^2 (2).
       \label{dilog_neg_two}
    \end{eqnarray}
    It follows from this result and the first equation in (\ref{eqn21_draft}) that
    \begin{eqnarray}
       \operatorname{Li}_2 \biggr{(} \frac{2}{3}\biggr{)} = d_2 + \frac{\pi^2}{6} + \frac{1}{2} \log^2(2) - \frac{1}{2} \log^2(3).
    \end{eqnarray}
    Moreover, using again formula (13) of Lewin's book \cite[p.283]{Lewin1981} with $x = 2$ along with the values for $\operatorname{Li}_2 (2)$ and $\operatorname{Li}_2 (-2)$ stated in equations (\ref{dilog_at2}) and (\ref{dilog_neg_two}), respectively, one arrives at
    \begin{eqnarray}
       \operatorname{Li}_2(4) = - 2d_2 + \frac{\pi^2}{6} - \log^2(2) - 2i \pi \log(2),
       \label{dilog_four}
    \end{eqnarray}
    Then with the value of $\operatorname{Li}_2(4)$ known in terms of the dilogarithm constant $d_2$, $x=4$ in formula (9) of \cite[p. 283]{Lewin1981} results in
    \begin{eqnarray}
       \operatorname{Li}_2 \biggr{(} \frac{4}{3} \biggr{)} = 2d_2 + \frac{\pi^2}{3} + \log^2(2) - \frac{1}{2} \log^2(3) + i \pi \log(3) - 2i \pi \log(2).
       \label{dilog_four_thirds}
    \end{eqnarray}
    Also, earlier we discussed a form of Landen's identity for the dilogarithm (\ref{landen_dilog}), which with $z=\frac{1}{2}$ leads to
    \begin{eqnarray}
       \operatorname{Li}_2 \biggr{(} \frac{1}{3} \biggr{)} = - d_2 - \frac{1}{2}\log^2 \biggr{(} \frac{3}{2} \biggr{)}.
       \label{dilog_third}
    \end{eqnarray}\\
    
    \hspace{-\parindent}It is clear from the examples considered above that computing the value of the dilogarithm constant $d_2 = \operatorname{Li}_2(-\frac{1}{2})$ is truly an important problem that is unfortunately still open. In light of how many values of the dilogarithm will immediately be known once the dilogarithm constant $d_2$ has been determined in the same spirit as $\operatorname{Li}_2(\frac{1}{2})$, $\operatorname{Li}_2(-1)$, and $\operatorname{Li}_2(1)$, we are hopeful the work done here will inspire other researchers to pursue this problem. 
    
\section{Bernoulli polynomials and solitons}
    \hspace{-\parindent}Let $B_p(x)\;$  ($p=0,1,2,\ldots)$ be the Bernoulli polynomials defined by the generating function
    \begin{eqnarray*}
        \frac{te^{xt}}{e^t - 1} = \sum_{n=0}^\infty B_n(x) \frac{t^n}{n!}
    \end{eqnarray*}
    ($B_p(0) = B_p$ are the Bernoulli numbers) \cite{Bateman1953, Berndt1985, Boyadzhiev2007, Sitaramachandrarao1987, Srivastava2001}.\\
    
    \hspace{-\parindent}A very interesting formula appears in section 7.12 of Lewin's book \cite{Lewin1981}, namely, equation (7.192) mentioned in the introduction
    \begin{eqnarray*}
        \operatorname{Li}_n(z) + (-1)^n \operatorname{Li}_n\biggr{(}\frac{1}{z}\biggr{)} = \frac{-2 \pi i}{n!} B_n \biggr{(} \frac{\log(z)}{2 \pi i} \biggr{)} \hspace{3mm} (n>1).
    \end{eqnarray*}
    The formula appeared previously also in \cite[1.11(8)]{Bateman1953}. The proof is based on a formula of Jonquiere which, on its part, is proved by contour integration (see \cite{Bateman1953}, \cite[(7.190)]{Lewin1981}, \cite{Pickard1968}). However, the formula in this form in both places \cite{Bateman1953, Lewin1981} is incorrect. For example, with $z=1$ and $n = 2p$ this equation gives
    \begin{eqnarray*}
        2 \zeta (2p) = \frac{-2\pi i}{(2p)!} B_{2p}
    \end{eqnarray*}
    while the correct formula of Euler for $\zeta(2p)$ is 
    \begin{eqnarray*}
        2 \zeta(2p) = \frac{(-1)^{p-1}(2\pi)^{2p}}{(2p)!} B_{2p}.
    \end{eqnarray*}
    Here we give a simple proof for the correct formula using the Fourier expansions of the Bernoulli polynomials.\\
    
    \hspace{-\parindent}\textbf{Proposition 3.} \textit{For every $p=1,2,\ldots$ we have the equations}
    
    \begin{eqnarray}
       \operatorname{Li}_{2p}(x) + \operatorname{Li}_{2p}(x^{-1}) = \frac{(-1)^{p+1} (2\pi)^{2p}}{(2p)!} B_{2p}\biggr{(} \frac{\log(x)}{2\pi i}\biggr{)}
       \label{sum_polylog_even}
    \end{eqnarray}
    \begin{eqnarray}
       \operatorname{Li}_{2p + 1} (x) - \operatorname{Li}_{2p+1} (x^{-1}) = \frac{(-1)^{p+1} (2\pi)^{2p+1}i}{(2p+1)!} B_{2p+1} \biggr{(}\frac{\log(x)}{2\pi i} \biggr{)}.
       \label{dif_polylog_odd}
    \end{eqnarray}
    
    \hspace{-\parindent}\textit{Proof.} Given $\operatorname{Li}_p(z) = \sum_{n=1}^\infty \frac{z^n}{n^p}$ we write
    \begin{eqnarray*}
       \operatorname{Li}_p(z) + \operatorname{Li}_p(z^{-1}) = \sum_{n=1}^\infty \frac{z^n + z^{-n}}{n^p} \hspace{3mm}\text{and} \hspace{3mm} \operatorname{Li}_p(z) - \operatorname{Li}_p(z^{-1}) = \sum_{n=1}^\infty \frac{z^n - z^{-n}}{n^p} 
    \end{eqnarray*}
    with $z = e^{2\pi  i t}$ we get the two equations
    \begin{eqnarray*}
        &&\operatorname{Li}_p(e^{2\pi i t}) + \operatorname{Li}_p(e^{-2 \pi i t}) = 2 \sum_{n=1}^\infty \frac{\cos(2 \pi n t)}{n^p} \\
        &&\operatorname{Li}_p(e^{2 \pi i t}) - \operatorname{Li}_p(e^{-2 \pi i t}) = 2 i \sum_{n=1}^\infty \frac{\sin(2 \pi n t)}{n^p}.
    \end{eqnarray*}
    At the same time we have the Fourier series expansions \cite[1.13]{Bateman1953}, \cite[p.65]{Srivastava2001}, \cite[p. 87]{Srivastava2011} for $0 \leq t \leq 1$
    \begin{eqnarray*}
        &&B_{2p+1}(t) = \frac{(-1)^{p+1}(2p+1)!}{2^{2p}\pi^{2p + 1}} \sum_{n=1}^\infty \frac{\sin(2 \pi n t)}{n^{2p+1}} \\
        &&B_{2p}(t) = \frac{(-1)^{p+1}(2p)!}{2^{2p-1} \pi^{2p}} \sum_{n=1}^\infty \frac{\cos(2 \pi n t)}{n^{2p}}.
    \end{eqnarray*}
    From here
    \begin{eqnarray*}
        &&2i \sum_{n=1}^\infty \frac{\sin(2 \pi n t)}{n^{2p+1}} = \frac{(-1)^{p+1} (2 \pi)^{2p+1}i}{(2p+1)!} B_{2p+1}(t)\\
        &&2 \sum_{n=1}^\infty \frac{\cos(2\pi n t)}{n^{2p}} = \frac{(-1)^{p+1} (2\pi)^{2p}}{(2p)!} B_{2p}(t) 
    \end{eqnarray*}
    (see also \cite[(17.4.4), (17.4.5)]{Hansen1975}). This gives
    \begin{eqnarray}
       &&\operatorname{Li}_{2p}(e^{2\pi i t}) + \operatorname{Li}_{2p}(e^{-2\pi i t}) = \frac{(-1)^{p+1} (2\pi)^{2p}}{(2p)!} B_{2p}(t) \label{eqn43_draft} \\
       &&\operatorname{Li}_{2p+1}(e^{2\pi i t}) - \operatorname{Li}_{2p+1}(e^{-2\pi i t}) = \frac{(-1)^{p+1} (2\pi)^{2p+1}i}{(2p+1)!} B_{2p+1}(t) \label{eqn44_draft}
    \end{eqnarray}
    and with the substitution $t = \frac{\log(x)}{2 \pi i}$ we come to the desired equations (\ref{sum_polylog_even}) and (\ref{dif_polylog_odd}). \\
    
    \hspace{-\parindent}Now consider the function $u(x,t) = -2 \sech^2 (x-4t)$. This function is the important soliton solution to the Korteweg-de Vries equation $u_t - 6 u u_x + u_{xxx} = 0$ \cite{Boyadzhiev2007}. It was shown in \cite{Boyadzhiev2007} that for every $n \geq 0$
    \begin{eqnarray}
       \int_{-\infty}^\infty x^n \sech^2(x-t) dx = 2(-i)^n \pi^n B_n \biggr{(}\frac{1}{2} + \frac{it}{\pi} \biggr{)}
       \label{integral_khristo2007}
    \end{eqnarray}
    expressing the moments of the soliton in terms of Bernoulli polynomials. Now in view of (\ref{sum_polylog_even}) we can express these moments in terms of polylogarithms. \\
    
    \hspace{-\parindent}\textbf{Corollary 4.}
    \begin{eqnarray}
       \int_{-\infty}^\infty x^{2p} \sech^2(x-t) dx = \frac{(-1)^{p+1}(2p)!}{2^{2p-1}} \Biggr{(} \operatorname{Li}_{2p}(-e^{-2it}) + \operatorname{Li}_{2p}(-e^{2it})\Biggr{)}
       \label{eqn46_draft}
    \end{eqnarray}
    \begin{eqnarray}
       \int_{-\infty}^\infty x^{2p+1} \sech^2(x-t) dx = \frac{i(2p+1)!}{2^{2p}} \Biggr{(} \operatorname{Li}_{2p+1}(-e^{-2it}) - \operatorname{Li}_{2p+1}(-e^{2it}) \Biggr{)}
       \label{eqn47_draft}
    \end{eqnarray}
    
    \hspace{-\parindent}\textit{Proof.} The substitution $t \to \frac{1}{2} + \frac{it}{\pi}$ brings (\ref{eqn43_draft}) to (\ref{eqn46_draft}) and (\ref{eqn44_draft}) to (\ref{eqn47_draft}) in view of (\ref{integral_khristo2007}).\\
    
    \hspace{-\parindent}With $t=0$ in (\ref{eqn47_draft}) we find $\int_{-\infty}^\infty x^{2p+1} \sech^2(x) dx = 0$ (as expected with odd integrand). With $t=0$ in (\ref{eqn46_draft}) we come to
    \begin{eqnarray*}
       \int_{-\infty}^\infty x^{2p} \sech^2(x) dx & = & \frac{(-1)^{p+1}(2p)!}{2^{2p-1}} \biggr{(}2 \operatorname{Li}_{2p}(-1) \biggr{)} \\
       & = & \frac{2(-1)^{p+1}(2p)!}{2^{2p-1}} \biggr{(}  \frac{1}{2^{2p-1}} - 1 \biggr{)} \zeta(2p). 
    \end{eqnarray*}
    
    \bibliography{main}
    \bibliographystyle{plain}

\end{document}